# OPTIMAL ADAPTIVE ESTIMATION OF A QUADRATIC FUNCTIONAL[1]

By T. Tony Cai and Mark G. Low

*University of Pennsylvania*

Adaptive estimation of a quadratic functional over both Besov and $L_p$ balls is considered. A collection of nonquadratic estimators are developed which have useful bias and variance properties over individual Besov and $L_p$ balls. An adaptive procedure is then constructed based on penalized maximization over this collection of nonquadratic estimators. This procedure is shown to be optimally rate adaptive over the entire range of Besov and $L_p$ balls in the sense that it attains certain constrained risk bounds.

**1. Introduction.** The problem of estimating the quadratic functional $\int f^2$ has received much attention in the statistical literature especially since, in a density estimation setting, Bickel and Ritov [5] showed under Hölder smoothness conditions that there is a breakdown in the minimax rate of convergence. Fully efficient estimation is possible when the function satisfies a Hölder smoothness condition with $\alpha > \frac{1}{4}$. However when $\alpha \le \frac{1}{4}$ minimax rates of convergence under mean squared error are of the order $n^{-8\alpha/(1+4\alpha)}$.

This theory has been developed and extended in a number of important directions which can be particularly easily described for the Gaussian sequence model

$$Y_i = \theta_i + n^{-1/2} z_i, \qquad i = 1, 2, \ldots, \tag{1}$$

where $z_i$ are i.i.d. standard normal random variables and where $\theta = (\theta_1, \theta_2, \ldots)$ is assumed to belong either to an $L_p$ or a Besov ball. Such models occupy a central role in the nonparametric function estimation literature. See, for example, [22]. In this sequence model setting estimation of the quadratic functional $Q(\theta) = \sum_{i=1}^{\infty} \theta_i^2$ is the analog of estimating the functional $\int f^2$ in the density estimation model.

Received November 2004; revised July 2005.

[1]Supported in part by NSF Grant DMS-03-06576.

The $L_p$ balls are defined as

$$L_p(\alpha, M) = \left\{\theta : \left(\sum_{i=1}^{\infty} i^{ps}|\theta_i|^p\right)^{1/p} \le M\right\}, \tag{2}$$

where $p > 0$, $\alpha > 0$, $M > 0$ and $s = \alpha + \frac{1}{2} - \frac{1}{p} > 0$. Besov balls in sequence space are typically defined in terms of a doubly indexed sequence $\{\theta_{j,k} : j = 0, 1, \ldots, k = 0, \ldots, 2^j - 1\}$. For $p, q, \alpha, M > 0$ the Besov ball $B^{\alpha}_{p,q}(M)$ is then given by

$$B^{\alpha}_{p,q}(M) = \left\{\theta : \left(\sum_{j=0}^{\infty}\left(2^{js}\left(\sum_{k=0}^{2^j-1}|\theta_{j,k}|^p\right)^{1/p}\right)^q\right)^{1/q} \le M\right\}, \tag{3}$$

where once again $s = \alpha + \frac{1}{2} - \frac{1}{p} > 0$. In particular, Besov balls contain as special cases a number of well-known smoothness spaces such as Hölder and Sobolev balls. It is possible to give a unified treatment of Besov balls and $L_p$ balls by setting in the case of Besov balls $\theta_i = \theta_{j,k}$, where $i = 2^j + k$. Noisy observation of Besov coefficients can then still be written as in (1). This convention is used throughout the paper, where in addition we shall assume that $p, q, \alpha, s > 0$.

For estimation of the quadratic functional $Q(\theta)$ over Besov and $L_p$ balls there are really two distinct cases of interest. The "dense" case corresponds to $p \ge 2$ and the "sparse" case to $p < 2$. Previous literature has focused primarily on the dense case where the parameter space is quadratically convex. In such cases the minimax theory for estimating the quadratic functional $Q(\theta)$ was well developed in [13] and [15]. In particular, this theory covers Besov balls $B^{\alpha}_{p,q}(M)$ and $L_p$ balls $L_p(\alpha, M)$ when $p \ge 2$. An important feature of this minimax theory is that optimal quadratic rules can be found within a "small" constant factor of the minimax risk.

The minimax theory for parameter spaces which are not quadratically convex is quite different. The near minimaxity of optimal quadratic rules typically does not hold when the parameter space is not quadratically convex. Cai and Low [10] develop the minimax theory in such cases over all Besov balls and $L_p$ balls with $p < 2$. A nonquadratic minimax procedure is given based on term-by-term thresholding. The nonquadratic procedure is sometimes fully efficient even when optimal quadratic rules have slow rates of convergence.

The minimax results for estimating the quadratic functional $Q(\theta)$ over $\Theta = L_p(\alpha, M)$ or $\Theta = B^{\alpha}_{p,q}(M)$ can be summarized as follows. Set $p_* = \min\{p, 2\}$, $s_* = \alpha + \frac{1}{2} - \frac{1}{p_*}$ and let

$$r(\alpha,p)=\begin{cases}1, & \text{if } \alpha p_* \ge \frac{1}{2},\\ 2-\dfrac{p_*}{1+2p_*s_*}, & \text{if } \alpha p_* < \frac{1}{2}.\end{cases} \tag{4}$$

Then

$$\inf_{\hat{Q}} \sup_{\theta\in\Theta} E_\theta(\hat{Q}-Q(\theta))^2 \asymp n^{-r(\alpha,p)}. \tag{5}$$

Moreover, if $\alpha p_* > \frac{1}{2}$,

$$\inf_{\hat{Q}} \sup_{\theta\in\Theta} E_\theta(\hat{Q}-Q(\theta))^2 = 4A(\Theta)n^{-1}(1+o(1)), \tag{6}$$

where $A(\Theta)=\sup_{\theta\in\Theta}\sum_{i=1}^{\infty}\theta_i^2$ is a constant and $4A(\Theta)n^{-1}$ is the inverse of the nonparametric Fisher information.

In comparison to minimax theory, the theory of adaptive estimation of $Q(\theta)$ is not as well developed. Most of the progress has been made in quadratically convex cases. Efromovich and Low [14] considered adaptive estimation of $Q(\theta)$ over hyper-rectangles, which corresponds to $L_p$ balls with $p=\infty$. It was shown that rate optimal adaptive estimators do not exist and that logarithmic penalties must be paid. An adaptive procedure only paying these logarithmic penalties was constructed using the method due to Lepski [26]. Tribouley [29] and Johnstone [21] developed an alternative adaptive procedure based on block thresholding algorithms. Gayraud and Tribouley [18] also used a block thresholding scheme for adaptation over Besov spaces with $p=2$ and $q=\infty$. Using Lepski's method to choose within a collection of quadratic rules Klemelä [23] considered sharp adaptation for $L_p$ balls with $p>2$.

All of the results mentioned so far focus on quadratically convex cases where $p\ge 2$. The sparse case where $p<2$ presents some major new difficulties which requires a novel approach for the construction of adaptive procedures. The goal of the present work is to develop a procedure which adapts simultaneously over all Besov and $L_p$ balls. This problem is significantly different from adaptation only over the dense cases where one can select from a collection of quadratic estimators. In the sparse case even minimax theory requires nonquadratic rules.

It is well known from previous work that block thresholding is an effective tool for adaptive estimation of $Q(\theta)$ in the dense case. Block thresholding can be used to guard against the worst case when there are a large number of small coefficients and where the exact location of these coefficients is unknown. On the other hand, in the sparse case, as shown in [10], the worst case occurs when there are a relatively small number of large coefficients with unknown location. In such cases term-by-term thresholding is effective. Unfortunately term-by-term thresholding does not work well for the dense case, and likewise, block thresholding does not work well in the sparse case. In order to develop a procedure that can adapt simultaneously over both the sparse and the dense cases we incorporate both approaches.

There are three parts to the adaptive estimator given in this paper. The initial component is based on a simple unbiased quadratic estimate of the first part of the quadratic functional. The third component is based on a term-by-term thresholding procedure with slowly growing threshold. The most important component, at least for the sparse cases, is an estimate of the middle part of the quadratic functional $Q(\theta)$. This estimate is based on penalized maximization over a collection of estimators, each of which uses both block thresholding and term-by-term thresholding. We show that the resulting procedure simultaneously attains the benchmarks for adaptive estimation given in Section 2. In particular, it is fully efficient over the largest collection of Besov and $L_p$ balls for which efficient estimators exist while paying minimal penalty over all other Besov and $L_p$ balls.

More precisely, it follows from the theorems in the paper that the adaptive estimator satisfies for $\Theta = L_p(\alpha, M)$ or $\Theta = B^{\alpha}_{p,q}(M)$

$$\sup_{\theta\in\Theta} E_\theta(\hat{Q} - Q(\theta))^2 \le Cn^{-r(\alpha,p)}(\log n)^{2p_* s_*/(1+2p_* s_*)} \tag{7}$$

when $\alpha p_* \le \frac{1}{2}$ and

$$\sup_{\theta\in\Theta} E_\theta(\hat{Q} - Q(\theta))^2 = 4A(\Theta)n^{-1}(1+o(1)) \tag{8}$$

when $\alpha p_* > \frac{1}{2}$, where $A(\Theta) = \sup_{\theta\in\Theta}\sum_{i=1}^{\infty}\theta_i^2$. In other words, the estimator is adaptively fully efficient over all Besov bodies where efficient estimation is possible and only pays a logarithmic penalty when the minimax rate is slower than $n^{-1}$. In fact, it is also shown in the present paper that the upper bound given in (7) is rate sharp.

It is interesting to compare these results with those of an estimator based on model selection given in [25]. Their procedure, say $\hat{Q}_{\mathrm{LM}}$, maximizes penalized quadratic estimators. It was shown that for $p < 2$,

$$\begin{aligned} &\sup_{\theta\in L_p(\alpha,M)} E_\theta(\hat{Q}_{\mathrm{LM}} - Q(\theta))^2 \\ &\qquad \le C\min\left\{\left(\frac{\log n}{n^2}\right)^{4s/(1+4s)}, \left(\frac{\log n}{n}\right)^{4\alpha/(1+2\alpha)}\right\} + \frac{C}{n} \end{aligned} \tag{9}$$

for some constant $C > 0$. A comparison shows that these upper bounds are always larger than the upper bounds given in the present paper whenever the estimator $\hat{Q}_{\mathrm{LM}}$ of [25] has a bound larger than $O(n^{-1})$. A more detailed comparison is given in Section 2.

The paper is organized as follows. In Section 2 we develop benchmarks for the evaluation of adaptive estimators. The major focus of the paper is the construction of an adaptive estimator which is described in detail in Section 3. A collection of nonquadratic estimators with specific bias variance properties is constructed. The adaptive procedure is then built by selecting a penalized estimator over this collection through maximization. We show that the procedure simultaneously attains the benchmarks given in Section 4 over all Besov and $L_p$ balls. In particular, it is fully efficient over the largest collection of Besov and $L_p$ balls for which efficient estimators exist while paying minimal penalty over all other Besov and $L_p$ balls. Proofs are given in Section 4.

**2. The cost of adaptation in the sparse case.** The primary goal of the present work is to construct estimators of $Q(\theta) = \sum_{i=1}^{\infty} \theta_i^2$ which are adaptive over all Besov and $L_p$ balls. This goal, however, needs to be made precise because even in the dense case of $p \geq 2$ it is well known that fully minimax rate optimal adaptation of the quadratic functional $Q(\theta) = \sum_{i=1}^{\infty} \theta_i^2$ is not possible. See, for example, [14]. A penalty must be made over $L_p$ or Besov balls with $p \geq 2$ and $\alpha < \frac{1}{4}$. Hence in the present context a rate adaptive estimator is one which attains well defined lower bounds.

In this section we shall develop the appropriate lower bounds needed as a benchmark for the evaluation of adaptive procedures which are given in Section 3.2. We shall see that an entirely similar phenomenon occurs in the sparse case, although the exponent of the logarithmic penalty is different. In particular, the following theorem shows that fully rate adaptive estimation of the quadratic functional $Q(\theta) = \sum_{i=1}^{\infty} \theta_i^2$ is not possible over any pair of Besov spaces which have different minimax rates of convergence. In the following theorem denote by 0 the zero vector. Then $E_0$ denotes the expectation under the sequence model (1) when $\theta = 0$.

THEOREM 1. *Let $\hat{Q}$ be an estimator of the quadratic functional $Q(\theta) = \sum_{i=1}^{\infty} \theta_i^2$. Let $r(\alpha, p)$ be the minimax rate for estimating $Q(\theta)$ over $\Theta = B_{p,q}^{\alpha}(M)$ or $\Theta = L_p(\alpha, M)$. Suppose that*

$$E_0(\hat{Q} - Q(0))^2 \leq Cn^{-\gamma} \tag{10}$$

*for some constants $\gamma > r(\alpha, p)$ and $C > 0$. Then the maximum squared bias over $\Theta$ satisfies, for some constant $C' > 0$,*

$$\sup_{\theta \in \Theta}(E_\theta \hat{Q} - Q(\theta))^2 \geq C' n^{-r(\alpha,p)} (\log n)^{2p_* s_*/(1+2p_* s_*)}. \tag{11}$$

The theorem makes clear that rate optimal estimators over one Besov or $L_p$ ball must pay a logarithmic penalty for the maximum risk over all Besov and $L_p$ balls which have slower minimax rates of convergence. In fact, as shown in (11), this logarithmic penalty must be paid in terms of maximum squared bias.

The major use of the lower bound given in the above theorem is as a benchmark for the development of an adaptive estimator. Adaptive estimators which attain these bounds must over each parameter space inflate the maximum bias over that parameter space. In the next section we shall use this fact to guide us in the development of estimators which are adaptive in the sense that they attain the lower bound given in Theorem 1.

The proof of this theorem also immediately yields the following corollary which shows the "inflexibility" of minimax rate optimal estimators, at least in cases where the minimax rate is slower than $n^{-1}$. In particular, there does not exist an estimator which attains the exact minimax rate of convergence over any pair of Besov or $L_p$ balls which have different minimax rates of convergence.

COROLLARY 1. *Let $\hat{Q}$ be a minimax rate optimal estimator of the quadratic functional $Q(\theta) = \sum_{i=1}^{\infty} \theta_i^2$ over $\Theta = B_{p,q}^{\alpha}(M)$ or $\Theta = L_p(\alpha, M)$ where the minimax rate $r(\alpha, p) < 1$. That is,*

$$\sup_{\theta \in \Theta} E_\theta (\hat{Q} - Q(\theta))^2 \le D n^{-r(\alpha,p)} \tag{12}$$

*for some $D > 0$. Then*

$$E_0 (\hat{Q} - Q(0))^2 \ge D' n^{-r(\alpha,p)} \tag{13}$$

*for some $D' > 0$ and hence*

$$\sup_{\theta \in \Theta'} E_\theta (\hat{Q} - Q(\theta))^2 \ge D' n^{-r(\alpha,p)}, \tag{14}$$

*where $\Theta'$ is any Besov or $L_p$ ball.*

The benchmark given in Theorem 1 is useful for the evaluation of adaptive procedures over parameter spaces which have a minimax rate of convergence slower than $n^{-1}$. On the other hand, over Besov and $L_p$ balls with $\alpha p_* > \frac{1}{2}$, the minimax risk given in (6) is another useful benchmark. Estimators attaining (6) can be termed efficient since they attain a nonparametric information bound as given, for example, in [4]. See also [10].

**3. The construction of an adaptive procedure.** The major goal of the present paper is the construction of an estimator of $Q(\theta)$ which adapts over all Besov and $L_p$ balls. The development of such an adaptive estimator can perhaps best be understood by breaking this construction into two stages. In the first stage a collection of nonquadratic estimators is constructed using both block thresholding and term-by-term thresholding. These estimators have precise bias and variance properties. More specifically, for a given Besov or $L_p$ ball when the minimax rate is slower than $n^{-1}$, one of the estimators in the collection has maximum squared bias attaining the lower bound given in (11) and which has variance smaller than the minimax risk. On the other hand, when fully efficient estimation is possible, one of the estimators has negligible bias and the variance attains the minimax lower bound. The construction of these nonquadratic estimators is given in Section 3.1.

These nonquadratic estimators are then used to build an adaptive procedure. At this stage the adaptive estimator is created by maximizing penalized versions of these nonquadratic estimators where the penalty is chosen to be a logarithmic factor of the standard deviation of each of these estimators.

The general approach of model selection via penalization has been shown to be effective for a number of adaptive function estimation problems. See, for example, [1, 3, 6, 25]. In particular, a major advance in estimating the quadratic functional $Q(\theta)$ was made in [25], where it was shown that maximizing penalized quadratic estimators of $Q(\theta)$ can yield a procedure which is adaptive over certain Besov and $L_p$ balls. It is shown in Section 3.2 that the procedure based on maximizing the penalized nonquadratic estimators is adaptive over all Besov and $L_p$ bodies. A comparison with the estimator of Laurent and Massart [25] is given in Section 3.3.

3.1. *Nonquadratic estimators with specific bias and variance properties.* We start with the construction of a collection of estimators which have precise bias and variance properties. These estimators incorporate both block and term-by-term thresholding. It is known that block thresholding estimators can perform well for dense cases, that is, when $p \geq 2$ and that term-by-term thresholding estimators can be minimax rate optimal for sparse cases, that is, when $p < 2$. By combining block thresholding and term-by-term thresholding, estimators can be constructed which trade bias and variance in very useful ways for both the dense and sparse cases. More specifically, for a given Besov or $L_p$ ball we build an estimator that has inflated maximum squared bias and reduced variance and which in particular attains the adaptive rate of convergence for mean squared error.

It is useful to break the problem of estimating $Q(\theta)$ into three components as follows. Let $m_0 = \frac{n}{(\log n)^2}$ and $m_k = 2^k m_0$ for $k \geq 1$. Divide the indices $i$ beyond $m_0$ into blocks of increasing sizes so that the $k$th block is of size $m_k$.

Let $J$ be the largest integer satisfying $2^J \le n$ and set

$$\xi_0 = \sum_{i=1}^{m_0} \theta_i^2, \qquad \xi_{\mathrm{mid}} = \sum_{i=m_0+1}^{m_J} \theta_i^2 \quad \text{and} \quad \xi_{\mathrm{tail}} = \sum_{i=m_J+1}^{\infty} \theta_i^2. \tag{15}$$

Then clearly $Q(\theta) = \xi_0 + \xi_{\mathrm{mid}} + \xi_{\mathrm{tail}}$. We shall use different strategies for estimating the three components $\xi_0$, $\xi_{\mathrm{mid}}$ and $\xi_{\mathrm{tail}}$. Estimation of $\xi_{\mathrm{mid}}$ is the most involved and so we shall first describe estimators for $\xi_0$ and $\xi_{\mathrm{tail}}$.

The component $\xi_0$ is naturally estimated by the unbiased quadratic estimator

$$\hat{\xi}_0 = \sum_{i=1}^{m_0} \left( Y_i^2 - \frac{1}{n} \right). \tag{16}$$

Note that for $\Theta = B_{p,q}^{\alpha}(M)$ or $\Theta = L_p(\alpha, M)$

$$\sup_{\theta \in \Theta} E_\theta (\hat{\xi}_0 - \xi_0)^2 = \sup_{\theta \in \Theta} \left\{ \frac{4\xi_0}{n} + \frac{2m_0}{n^2} \right\} = \frac{4A(\Theta)}{n} (1 + o(1)), \tag{17}$$

where $A(\Theta) = \sup_{\theta \in \Theta} \sum_{i=1}^{\infty} \theta_i^2$. It is clear that this term is equal to the minimax risk when fully efficient estimation is possible and negligible whenever the minimax rate of convergence is slower than the parametric rate of $n^{-1}$.

The technique underlying the estimation of the tail component is similar to that used for minimax estimation of $Q(\theta)$ in the sparse case as given in [10]. First define $\gamma_i$ by

$$\gamma_i = 2 \left( \left\lceil \log_2 \frac{i}{m_J} \right\rceil + 1 \right), \qquad i \ge m_J + 1, \tag{18}$$

where $\lceil x \rceil$ denotes the smallest integer greater than or equal to $x$. That is, $\gamma_i = 2(j - J + 2)$ for $m_j + 1 \le i \le m_{j+1}$ and $j \ge J$. Then the tail component $\xi_{\mathrm{tail}}$ is estimated by a term-by-term thresholding estimator with slowly growing threshold

$$\hat{\xi}_{\mathrm{tail}} = \sum_{i=m_J+1}^{\infty} \left( Y_i^2 - \frac{\gamma_i \log n}{n} \right)_+. \tag{19}$$

We shall show that the risk due to estimation of the tail is always negligible relative to the minimax risk for $\Theta = B_{p,q}^{\alpha}(M)$ or $\Theta = L_p(\alpha, M)$, that is,

$$\sup_{\theta \in \Theta} E_\theta (\hat{\xi}_{\mathrm{tail}} - \xi_{\mathrm{tail}})^2 = o(n^{-r(\alpha,p)}). \tag{20}$$

We now turn to estimation of the middle component $\xi_{\mathrm{mid}}$, which is more involved and uses both block thresholding and term-by-term thresholding. Let $\Theta = B_{p,q}^{\alpha}(M)$ or $\Theta = L_p(\alpha, M)$. The estimator $\hat{\xi}_{\mathrm{mid}}$ depends on the parameters $\alpha$ and $p$. For each integer $k$ such that $1 \le k \le J-1$ set $\tau_{k,i} = 2(j+1-k)$ for $m_j + 1 \le i \le m_{j+1}$ and $k \le j \le J-1$. That is,

$$\tau_{k,i} = 2\left\lceil \log_2 \frac{i}{m_k} \right\rceil, \qquad i \ge m_k + 1, \tag{21}$$

where $\lceil x \rceil$ once again denotes the smallest integer greater than or equal to $x$. For $i \ge m_k + 1$, set $\mu_{k,i} = E_0\{(Y_i^2 - \frac{\tau_{k,i}}{n})_+\}$ where the expectation is taken under $\theta = 0$. Let

$$\lambda_k = \frac{(m_k - m_0) + 2\sqrt{(m_k - m_0)\log(m_k - m_0)}}{n}.$$

For each $1 \le k \le J-1$ set

$$\hat{\xi}_k = \left(\sum_{i=m_0+1}^{m_k} Y_i^2 - \lambda_k\right)_+ + \sum_{i=m_k+1}^{m_J} \left[\left(Y_i^2 - \frac{\tau_{k,i}}{n}\right)_+ - \mu_{k,i}\right]. \tag{22}$$

Recall that $p_* = \min\{p, 2\}$ and $s_* = \alpha + \frac{1}{2} - \frac{1}{p_*}$. Set $k_*$ to be the largest integer such that

$$m_{k_*} = 2^{k_*} m_0 \le \max\{2m_0, n^{p_*/(1+2p_* s_*)} (\log n)^{-1/(1+2p_* s_*)}\}. \tag{23}$$

The middle component $\xi_{\text{mid}}$ is then estimated by

$$\begin{aligned} \hat{\xi}_{\text{mid}} &= \hat{\xi}_{k_*} \\ &= \left(\sum_{i=m_0+1}^{m_{k_*}} Y_i^2 - \lambda_{k_*}\right)_+ + \sum_{i=m_{k_*}+1}^{m_J} \left\{\left(Y_i^2 - \frac{\tau_{k_*,i}}{n}\right)_+ - \mu_{k_*,i}\right\}. \end{aligned} \tag{24}$$

We shall show that the risk of $\hat{\xi}_{\text{mid}}$ for estimating $\xi_{\text{mid}}$ is negligible when fully efficient estimation is possible and otherwise attains the lower bound given in Theorem 1 over the given $\Theta$. It should also be noted that the first term used to define $\hat{\xi}_{\text{mid}}$ would suffice for the dense case where $p \ge 2$. The second term is needed for the sparse case where $p < 2$.

The quadratic functional $Q(\theta)$ is then estimated by

$$\hat{Q}_{k_*} = \hat{\xi}_0 + \hat{\xi}_{\text{mid}} + \hat{\xi}_{\text{tail}}. \tag{25}$$

The following result shows that this estimator has desirable bias and variance properties.

PROPOSITION 1. *Let $\Theta = B_{p,q}^{\alpha}(M)$ or $\Theta = L_p(\alpha, M)$ and let the estimator $\hat{Q}_{k_*}$ be given as in (25). If $\alpha p_* < \frac{1}{2}$, then the maximum squared bias satisfies*

$$\sup_{\theta \in \Theta} (E\hat{Q}_{k_*} - Q(\theta))^2 \le C n^{-(2 - p_*/(1+2p_* s_*))} (\log n)^{2p_* s_*/(1+2p_* s_*)} \tag{26}$$

*and the maximum variance satisfies*

$$\sup_{\theta\in\Theta} \mathrm{Var}(\hat{Q}_{k_*}) \le C n^{-(2-p_*/(1+2p_*s_*))} (\log n)^{-1/(1+2p_*s_*)}. \tag{27}$$

*On the other hand, if* $\alpha p_* > \frac{1}{2}$, *then* $\hat{Q}_{k_*}$ *is asymptotically efficient, that is,*

$$\sup_{\theta\in\Theta} E_\theta(\hat{Q}_{k_*} - Q(\theta))^2 = 4A(\Theta)\, n^{-1}(1+o(1)), \tag{28}$$

*where* $A(\Theta) = \sup_{\theta\in\Theta}\sum_{i=1}^{\infty}\theta_i^2$. *Furthermore, in the boundary case of* $\alpha p_* = \frac{1}{2}$, $\hat{Q}_{k_*}$ *satisfies*

$$\sup_{\theta\in\Theta} (E\hat{Q}_{k_*} - Q(\theta))^2 \le C n^{-1} (\log n)^{2p_*s_*/(1+2p_*s_*)} \tag{29}$$

*and*

$$\sup_{\theta\in\Theta} \mathrm{Var}(\hat{Q}_{k_*}) = 4A(\Theta)\, n^{-1}(1+o(1)), \tag{30}$$

*where once again* $A(\theta) = \sup_{\theta\in\Theta}\sum_{i=1}^{\infty}\theta_i^2$.

Note that the estimator $\hat{Q}_{k_*}$ has reduced variance and inflated bias compared to the minimax risk when the minimax rate of convergence is slower than the parametric rate. In fact in these cases the ratio of the maximum squared bias to the maximum variance is exactly of order $\log n$. These properties are crucial in the construction of the adaptive estimator given in Section 3.2.

3.2. *Adaptive procedure.* We shall now turn to the construction of a general adaptive procedure building upon the collection of nonquadratic estimators given in Section 3.1. The adaptive estimator is the maximization of penalized versions of these nonquadratic estimators. Let $\hat{Q}_k = \hat{\xi}_0 + \hat{\xi}_k + \hat{\xi}_{\mathrm{tail}}$ where $\hat{\xi}_0$, $\hat{\xi}_k$, and $\hat{\xi}_{\mathrm{tail}}$ are defined in (16), (22) and (19), respectively. The adaptive estimator is then given by

$$\hat{Q} = \max_{1\le k\le J}\left\{\hat{Q}_k - \frac{6\sqrt{m_k \log n}}{n}\right\}. \tag{31}$$

We shall show later that for any given Besov or $L_p$ ball the penalty term in (31) is always a logarithmic factor larger than the maximum variance of the estimator $\hat{Q}_k$. Moreover, the bias of the estimators $\hat{Q}_k$ is always negligible when it is positive, whereas in worst cases it must be negative. Taking a maximization with the penalty term results in an optimal trading of bias and variance over all Besov and $L_p$ balls.

It is also convenient to define

$$\hat{\xi}_{\mathrm{mid}} = \max_{1\le k\le J}\left\{\hat{\xi}_k - \frac{6\sqrt{m_k \log n}}{n}\right\}. \tag{32}$$

Then the estimator $\hat{Q}$ can be equivalently written as

$$\hat{Q} = \hat{\xi}_0 + \hat{\xi}_{\mathrm{mid}} + \hat{\xi}_{\mathrm{tail}}, \tag{33}$$

where $\hat{\xi}_0$, $\hat{\xi}_{\mathrm{mid}}$ and $\hat{\xi}_{\mathrm{tail}}$ are defined in (16), (32) and (19), respectively.

The following theorem shows that the estimator $\hat{Q}$ is optimally adaptive over all Besov and $L_p$ balls both for the dense and sparse cases.

THEOREM 2. *Let $Q(\theta) = \sum_{i=1}^{\infty} \theta_i^2$ and let the estimator $\hat{Q}$ be defined as in (31). Then the risk of $\hat{Q}$ satisfies for all Besov balls $\Theta = B_{p,q}^{\alpha}(M)$ and all $L_p$ balls $\Theta = L_p(\alpha, M)$*

$$\begin{aligned} &\sup_{\theta\in\Theta} E_\theta(\hat{Q} - Q(\theta))^2 \\ &\quad \le \begin{cases} 4A(\Theta)\, n^{-1}(1 + o(1)), & \text{for } \alpha p_* > \frac{1}{2}, \\ Cn^{-(2-p_*/(1+2p_*s_*))}(\log n)^{2p_*s_*/(1+2p_*s_*)}, & \text{for } \alpha p_* \le \frac{1}{2}, \end{cases} \end{aligned} \tag{34}$$

*where $C > 0$ and $A(\Theta) = \sup_{\theta\in\Theta} \sum_{i=1}^{\infty} \theta_i^2$ are constants.*

Comparing the upper bounds given in the above theorem with the lower bound given in Theorem 1 as well as the information bound given in (6), it is clear that the estimator $\hat{Q}$ is adaptive over all Besov and $L_p$ balls. In particular, it is adaptively efficient over those parameter spaces where efficient estimation is possible.

3.3. *Discussion.* It is interesting to compare the performance of the adaptive estimator $\hat{Q}$ with the estimator, say $\hat{Q}_{\mathrm{LM}}$, given in [25]. The estimator there is constructed based on model selection. It chooses a penalized quadratic estimator through maximization. In contrast, in this paper the adaptive estimator $\hat{Q}$ selects among a collection of penalized nonquadratic estimators. These nonquadratic estimators enable optimal adaptation over sparse cases corresponding to Besov and $L_p$ balls with $p < 2$ in addition to the standard dense cases of $p \ge 2$. In Table 1, $R(\hat{Q}_{\mathrm{LM}})$ denotes the order of the risk upper bound of $\hat{Q}_{\mathrm{LM}}$ given in [25] and $R(\hat{Q})$ denotes the order of the maximum risk of $\hat{Q}$ as given in Theorem 2. The comparison is focused on the sparse case where $p < 2$.

Simple algebra shows that the risk upper bounds for $\hat{Q}_{\mathrm{LM}}$ are always larger by an algebraic factor than those for $\hat{Q}$ whenever $R(\hat{Q}_{\mathrm{LM}}) \gg n^{-1}$. In particular, if $1 \le p < \frac{4}{3}$ and $\frac{1}{2p} < \alpha \le \frac{1}{2}$, $R(\hat{Q}_{\mathrm{LM}})$ is of order $(\frac{\log n}{n})^{4\alpha/(1+2\alpha)}$ whereas $\hat{Q}$ is fully efficient. Likewise when $\frac{4}{3} \le p < 2$ and $\frac{1}{2p} < \alpha \le \frac{1}{p} - \frac{1}{4}$, $R(\hat{Q}_{\mathrm{LM}})$ is of order $(\frac{\log n}{n^2})^{4s/(1+4s)}$ and $\hat{Q}$ is once again fully efficient.

It is also interesting to note that the problem of estimating the quadratic functional $Q(\theta)$ is strongly connected to the problem of estimating linear functionals. This connection was developed first in [13] where it was shown that a modulus of continuity for orthosymmetric parameter spaces could be used to yield optimal quadratic minimax estimators in a way that is analogous to a similar theory for minimax estimation of linear functionals given in [12]. See [10] for further discussion of this connection and the connection to estimating the whole signal $\theta$ in the minimax estimation setting. The adaptation theory for estimating the quadratic functional $Q(\theta)$ developed in the present paper is also similar to that for estimating linear functionals. For linear functionals Lepski [26] was the first to show that logarithmic penalties must often be paid when adapting over collections of parameter spaces, as is the case in the present paper. Further refinements and generalizations for adaptive estimation of linear functionals can be found in [9, 24, 27].

**4. Proofs.** The main results are proved in the order of Proposition 1, Theorem 2 and then Theorem 1. Detailed proofs are only given for $L_p$ balls since the proofs for Besov balls are entirely analogous. In this section $C$ denotes a positive constant not depending on $n$ that may vary from place to place, $\phi(z)$ and $\Phi(z)$ denote the density and cumulative distribution function of a standard normal random variable and $\tilde{\Phi}(z) = 1 - \Phi(z)$.

4.1. *Preparatory results.* The following lemma helps in the analysis of term-by-term thresholding estimators and is important to the proof of Proposition 1 and Theorem 2.

LEMMA 1. *Let $X \sim N(\theta, \frac{1}{n})$ and $\tau \geq 1$. Set $\mu_0(\tau) = E_0\{(X^2 - \frac{\tau}{n})_+\}$ where the expectation is taken under $\theta = 0$. Let $\hat{\xi} = (X^2 - \frac{\tau}{n})_+ - \mu_0(\tau)$. Then*

TABLE 1
*Comparison of the performance of the estimators $\hat{Q}_{\mathrm{LM}}$ and $\hat{Q}$*

| | $\mathbf{0 < p < 1}$ | $\mathbf{1 \leq p < \frac{4}{3}}$ | | |
|---|---|---|---|---|
| | $\alpha > \frac{1}{p} - \frac{1}{2}$ | $\alpha \leq \frac{1}{2p}$ | $\frac{1}{2p} < \alpha \leq \frac{1}{2}$ | $\alpha > \frac{1}{2}$ |
| $R(\hat{Q}_{\mathrm{LM}})$ | $n^{-1}$ | $(\frac{\log n}{n})^{4\alpha/(1+2\alpha)}$ | $(\frac{\log n}{n})^{4\alpha/(1+2\alpha)}$ | $n^{-1}$ |
| $R(\hat{Q})$ | $n^{-1}$ | $n^{-(2-p/(1+2ps))} \times (\log n)^{2ps/(1+2ps)}$ | $n^{-1}$ | $n^{-1}$ |
| | | $\mathbf{\frac{4}{3} \leq p < 2}$ | | |
| | $\alpha \leq \frac{2}{p} - 1$ | $\frac{2}{p} - 1 < \alpha \leq \frac{1}{2p}$ | $\frac{1}{2p} < \alpha \leq \frac{1}{p} - \frac{1}{4}$ | $\alpha > \frac{1}{p} - \frac{1}{4}$ |
| $R(\hat{Q}_{\mathrm{LM}})$ | $(\frac{\log n}{n})^{4\alpha/(1+2\alpha)}$ | $(\frac{\log n}{n^2})^{4s/(1+4s)}$ | $(\frac{\log n}{n^2})^{4s/(1+4s)}$ | $n^{-1}$ |
| $R(\hat{Q})$ | $n^{-(2-p/(1+2ps))} \times (\log n)^{2ps/(1+2ps)}$ | $n^{-(2-p/(1+2ps))} \times (\log n)^{2ps/(1+2ps)}$ | $n^{-1}$ | $n^{-1}$ |

$0 < \mu_0(\tau) \le \frac{4}{\sqrt{2\pi} n \tau^{1/2} e^{\tau/2}}$,

$$|E_\theta \hat{\xi} - \theta^2| \le \min\left(\frac{2\tau}{n}, \theta^2\right) \tag{35}$$

*and the variance of* $\hat{\xi}$ *satisfies*

$$\mathrm{Var}(\hat{\xi}) \le \frac{6\theta^2}{n} + \frac{4\tau^{1/2} + 18}{n^2 e^{\tau/2}}. \tag{36}$$

*In addition, if* $Z \sim N(0,1)$ *and* $V(\tau) = \mathrm{Var}[(Z^2 - \tau)_+]$ *then*

$$V(\tau) \le (16\tau^{-1/2} - 9\tau^{-3/2} + 9\tau^{-5/2})\phi(\tau^{1/2}). \tag{37}$$

PROOF. Equations (35) and (36) are from [10]. For (37), it follows from the standard alternate series tail bound $\tilde{\Phi}(z) \le (z^{-1} - z^{-3} + 3z^{-5})\phi(z)$ for $z > 0$ that

$$\begin{aligned} V(\tau) &\le 2\int_{\tau^{1/2}}^{\infty} (z^2 - \tau)^2 \phi(z)\, dz \\ &= (6\tau^{1/2} - 2\tau^{3/2})\phi(\tau^{1/2}) + (6 - 4\tau + 2\tau^2)\tilde{\Phi}(\tau^{1/2}) \\ &\le (16\tau^{-1/2} - 9\tau^{-3/2} + 9\tau^{-5/2})\phi(\tau^{1/2}). \qquad \square \end{aligned}$$

Lemmas 2 and 3 given below provide useful properties of term-by-term thresholding estimators and are central to the proof of Theorem 2.

LEMMA 2. *Let* $X_i = \theta_i + Z_i$ *where* $Z_i \overset{i.i.d.}{\sim} N(0, \sigma^2)$ *for* $i = 1, 2, \ldots, m$. *Let* $\xi = \sum_{i=1}^m \theta_i^2$. *Let* $\lambda \ge 0$. *Then for any* $x$

$$\begin{aligned} &P\left(\sum_{i=1}^m (X_i^2 - \lambda)_+ \ge x\right) \\ &\qquad \le P\left\{[(Z_1 + \xi^{1/2})^2 - \lambda]_+ + \sum_{i=2}^m (Z_i^2 - \lambda)_+ \ge x\right\}. \end{aligned} \tag{38}$$

*That is, for a given value of* $\xi = \sum_{i=1}^m \theta_i^2$, *the random variables* $\sum_{i=1}^m (X_i^2 - \lambda)_+$ *are stochastically maximized when* $\theta_1 = \xi^{1/2}$ *and* $\theta_i = 0$ *for all* $i = 2, \ldots, m$.

PROOF. Intuitively the result of this lemma seems clear since given the sum $\sum_{i=1}^m X_i^2$ the value of $\sum_{i=1}^m (X_i^2 - \lambda)_+$ is a decreasing function of the number of nonzero terms in this sum. A formal proof can be given as follows. Begin with the case when $m = 2$. Let $x > 0$ and note that

$$\begin{aligned} &P\{(X_1^2 - \lambda)_+ + (X_2^2 - \lambda)_+ \ge x\} \\ &\qquad = E(P\{(X_1^2 - \lambda)_+ + (X_2^2 - \lambda)_+ \ge x | X_1^2 + X_2^2\}). \end{aligned}$$

It thus suffices to show that the conditional probability

$$g(x;\theta_1,\theta_2,\rho) = P\{(X_1^2-\lambda)_+ + (X_2^2-\lambda)_+ \geq x | X_1^2 + X_2^2 = \rho^2\}$$

is maximized when $\theta_1 = \xi^{1/2}$ and $\theta_2 = 0$ since the distribution of $X_1^2 + X_2^2$ depends on $\theta_1$ and $\theta_2$ only through $\theta_1^2 + \theta_2^2 = \xi$.

Note that if $\rho^2 \leq \lambda + x$, then $g(x;\theta_1,\theta_2,\rho) = 0 = g(x;\xi^{1/2},0,\rho)$. On the other hand, if $\rho^2 > 2\lambda + x$, then $g(x;\theta_1,\theta_2,\rho) = 1 = g(x;\xi^{1/2},0,\rho)$. Now consider the main case $\lambda + x < \rho^2 \leq 2\lambda + x$. In this case

$$\begin{aligned} g(x;\theta_1,\theta_2,\rho) = {} & P\{X_1^2 \geq \lambda + x | X_1^2 + X_2^2 = \rho^2\} \\ & + P\{X_2^2 \geq \lambda + x | X_1^2 + X_2^2 = \rho^2\}. \end{aligned}$$

It is more convenient to use polar coordinates by setting $X_1 = \rho\cos(\phi)$, $X_2 = \rho\sin(\phi)$. The conditional distribution of $(X_1, X_2)$ given $X_1^2 + X_2^2$ is a von Mises distribution. See, for example, [7, 28, 30].

Since the distribution of $X_i^2$ depends only on $\theta_i^2$, $i = 1,2$, without loss of generality we assume $\theta_i \geq 0$. Let $\beta$ be the angle between the direction of $(\theta_1,\theta_2)$ and the horizontal axis. More precisely, $\cos(\beta) = \frac{\theta_1}{\sqrt{\theta_1^2+\theta_2^2}}$. Then $0 \leq \beta \leq \frac{\pi}{2}$. The conditional distribution of $\phi$ given $\rho$ is thus given by $q_\beta(\phi) = ce^{d\cos(\phi-\beta)}$ where $c$ and $d$ are some positive constants. See [30]. Let

$$\begin{aligned} u_\beta(\phi) &= q_\beta(\phi) + q_\beta\left(\phi + \frac{\pi}{2}\right) + q_\beta(\phi+\pi) + q_\beta\left(\phi + \frac{3\pi}{2}\right) \\ &= ce^{d\cos(\phi-\beta)} + ce^{-d\cos(\phi-\beta)} + ce^{d\sin(\phi-\beta)} + ce^{-d\sin(\phi-\beta)}. \end{aligned}$$

Then $g(x;\theta_1,\theta_2,\rho) = 1 - \int_{-\phi_0}^{\phi_0} u_\beta(\phi)\,d\phi$, where $\phi_0 = \cos^{-1}\frac{\sqrt{\lambda+x}}{\rho}$. Note that $0 < \phi_0 \leq \frac{\pi}{4}$.

It is easy to check that $u_\beta(\phi)$ has the following properties:

- It is periodic with period $\frac{\pi}{2}$, $u_\beta(\phi) = u_\beta(\phi + \frac{\pi}{2})$.
- $u_\beta(\phi)$ attains its maximum when $\phi = \beta$.
- $u_\beta(\beta + x) = u_\beta(\beta - x)$.
- $u_\beta$ is decreasing on the interval $[\beta, \beta + \frac{\pi}{4})$.

Noting the properties of $u_\beta$, it now follows from the rearrangement result in [20], page 278, that the above integral is maximized when $\beta = 0$, which corresponds to $\theta_1 = \xi^{1/2}$ and $\theta_2 = 0$. Hence $g(x;\theta_1,\theta_2,\rho) \leq g(x;\xi^{1/2},0,\rho)$. This completes the proof for $m = 2$. The general case now follows by first conditioning on $X_1^2 + X_2^2, X_3, \ldots, X_m$ and then by induction. □

LEMMA 3. *Let* $Z_i \overset{i.i.d.}{\sim} N(0,1)$, $i = 1,\ldots,m_n$, *with* $m_n \geq n$. *Let* $\gamma > 0$ *be fixed. Let* $\tau_{n,i} \geq 0$, $\mu_{n,i} = E[(Z_i^2 - \tau_{n,i})_+]$, $\sigma_{n,i}^2 = \mathrm{Var}[(Z_i^2 - \tau_{n,i})_+]$ *and* $V_n = \sum_{i=1}^{m_n} \sigma_{n,i}^2$. *Then there exists some absolute constant* $c_* > 0$ *such that for all sufficiently large* $n$

$$E\left\{\left[\sum_{i=1}^{m_n}((Z_i^2-\tau_{n,i})_+-\mu_{n,i})-(\gamma V_n\log n)^{1/2}\right]_+\right\}^2 \tag{39}$$
$$\leq (2\gamma V_n \log n + c_* V_n^{1/2})(\gamma \log n)^{-1/4} n^{-\gamma/4}.$$

PROOF. Set $A_n = E\{(\sum_{i=1}^{m_n}[(Z_i^2-\tau_{n,i})_+-\mu_{n,i}]-(\gamma V_n\log n)^{1/2})_+\}^2$. The Cauchy–Schwarz inequality then yields

$$A_n \leq \left\{E\left(\sum_{i=1}^{m_n}[(Z_i^2-\tau_{n,i})_+-\mu_{n,i}]-(\gamma V_n\log n)^{1/2}\right)^4\right\}^{1/2}$$
$$\times\left\{P\left(\sum_{i=1}^{m_n}[(Z_i^2-\tau_{n,i})_+-\mu_{n,i}]>(\gamma V_n\log n)^{1/2}\right)\right\}^{1/2}.$$

It is easy to verify by direct calculations that

$$\sup_{\tau_{n,i}\geq 0}\frac{E|(Z_i^2-\tau_{n,i})_+-\mu_{n,i}|^3}{\sigma_{n,i}^2}<\infty \tag{40}$$

and that the characteristic functions of $(Z_i^2-\tau_{n,i})_+$ are analytic. It then follows from [17], page 553, and the standard normal tail bound $\tilde{\Phi}(z) \leq z^{-1}\phi(z)$ for $z > 0$ that there exists some $n_* > 0$ such that for all $n \geq n_*$

$$P\left(\sum_{i=1}^{m_n}[(Z_i^2-\tau_{n,i})_+-\mu_{n,i}]>(\gamma V_n\log n)^{1/2}\right)\leq(\gamma\log n)^{-1/2}n^{-\gamma/2}. \tag{41}$$

Set $B_n = E(\sum_{i=1}^{m_n}[(Z_i^2-\tau_{n,i})_+-\mu_{n,i}]-(\gamma V_n\log n)^{1/2})^4$. It then follows from Rosenthal's inequality ([19], page 23) that for some absolute constant $c_1 > 0$

$$B_n \leq 4E\left(\sum_{i=1}^{m_n}[(Z_i^2-\tau_{n,i})_+-\mu_{n,i}]\right)^4+4(\gamma V_n\log n)^2$$
$$\leq c_1\left\{\sum_{i=1}^{m_n}E[(Z_i^2-\tau_{n,i})_+-\mu_{n,i}]^4+V_n^2\right\}+4(\gamma V_n\log n)^2.$$

It is also easy to verify by direct calculations that

$$\sup_{\tau_{n,i}\geq 0}\frac{E|(Z_i^2-\tau_{n,i})_+-\mu_{n,i}|^4}{\sigma_{n,i}^2}<\infty \tag{42}$$

and hence $B_n \leq (4\gamma^2\log^2 n + c_1)V_n^2 + c_2 V_n$ for some absolute constant $c_2 > 0$. Therefore for some absolute constant $c_* > 0$, $A_n \leq (2\gamma V_n \log n + c_* V_n^{1/2}) \times (\gamma\log n)^{-1/4} n^{-\gamma/4}$ for all sufficiently large $n$. □

4.2. *Proof of Proposition 1.* The proof of Proposition 1 relies heavily on Lemma 1. Denote by $B(\theta)$ and $V(\theta)$ the bias and variance of $\hat{Q}_{k_*}$, respectively. Set $\Theta = L_p(\alpha, M)$. Let $\xi_0$, $\xi_{\mathrm{mid}}$ and $\xi_{\mathrm{tail}}$ be given as in (15). Set $\xi_{\mathrm{mid1}} = \sum_{i=m_0+1}^{m_{k_*}} \theta_i^2, \xi_{\mathrm{mid2}} = \sum_{i=m_{k_*}+1}^{m_J} \theta_i^2, \hat{\xi}_{\mathrm{mid1}} = (\sum_{i=m_0+1}^{m_{k_*}} Y_i^2 - \lambda_{k_*})_+$ and $\hat{\xi}_{\mathrm{mid2}} = \sum_{i=m_{k_*}+1}^{m_J} \{(Y_i^2 - \frac{\tau_{k_*,i}}{n})_+ - \mu_{k_*,i}\}$.

We shall consider the bias and variance separately. First consider the variance. Note that $m_0 = \frac{n}{(\log n)^2}$ and $\xi = \sum_{i=1}^{\infty} \theta_i^2 \le A(\Theta)$. Hence

$$\mathrm{Var}(\hat{\xi}_0) = \sum_{i=1}^{m_0} \mathrm{Var}(Y_i^2) = \frac{4\xi_0}{n} + \frac{2m_0}{n^2} \le \frac{4A(\Theta)}{n}(1 + o(1)). \tag{43}$$

Note that for any random variable $X$, $\mathrm{Var}((X)_+) \le \mathrm{Var}(X)$. See, for example, [10]. Hence $\mathrm{Var}(\hat{\xi}_{\mathrm{mid1}}) \le \sum_{i=m_0+1}^{m_{k_*}} \mathrm{Var}(Y_i^2) = \frac{4\xi_{\mathrm{mid1}}}{n} + \frac{2(m_{k_*} - m_0)}{n^2}$. Lemma 1 yields that $\mathrm{Var}(\hat{\xi}_{\mathrm{mid2}}) \le \frac{6\xi_{\mathrm{mid2}}}{n} + \sum_{i=m_{k_*}+1}^{m_J} \frac{4\tau_{k_*,i}^{1/2} + 18}{n^2 e^{(\tau_{k_*,i})/2}}$ and consequently

$$\begin{aligned}\mathrm{Var}(\hat{\xi}_{\mathrm{mid}}) &= \mathrm{Var}(\hat{\xi}_{\mathrm{mid1}}) + \mathrm{Var}(\hat{\xi}_{\mathrm{mid2}}) \\ &\le \frac{6\xi_{\mathrm{mid}}}{n} + \frac{2m_{k_*}}{n^2} + \sum_{i=m_{k_*}+1}^{m_J} \frac{4\tau_{k_*,i}^{1/2} + 18}{n^2 e^{(\tau_{k_*,i})/2}}.\end{aligned}$$

Note that for $m_{k_*+j-1} + 1 \le i \le m_{k_*+j}$ and $j \ge 1$, $\tau_{k_*,i} = 2j$. Hence

$$\begin{aligned}\sum_{i=m_{k_*}+1}^{m_J} \frac{4\tau_{k_*,i}^{1/2} + 18}{n^2 e^{(\tau_{k_*,i})/2}} &= \sum_{j=1}^{J-k_*} \sum_{i=m_{k_*+j-1}+1}^{m_{k_*+j}} (4(2j)^{1/2} + 18) e^{-j} n^{-2} \\ &= \sum_{j=1}^{J-k_*} (4(2j)^{1/2} + 18) e^{-j} 2^{j-1} m_{k_*} n^{-2} \\ &\le C m_{k_*} n^{-2}\end{aligned} \tag{44}$$

for some constant $C > 0$, since $\sum_{j=1}^{\infty} (4(2j)^{1/2} + 18) e^{-j} 2^{j-1} < \infty$. Therefore

$$\mathrm{Var}(\hat{\xi}_{\mathrm{mid}}) \le 6\xi_{\mathrm{mid}} n^{-1} + C m_{k_*} n^{-2}. \tag{45}$$

For the tail component, note that Lemma 1 once again yields

$$\mathrm{Var}(\hat{\xi}_{\mathrm{tail}}) \le \frac{6\xi_{\mathrm{tail}}}{n} + \sum_{i=m_J+1}^{\infty} \frac{4\gamma_i^{1/2} (\log n)^{1/2} + 18}{n^{2+\gamma_i/2}}.$$

Using similar derivation as in (44), we have $\sum_{i=m_J+1}^{\infty} \frac{4\gamma_i^{1/2} (\log n)^{1/2} + 18}{n^{2+\gamma_i/2}} \le Cn^{-2} \times (\log n)^{1/2}$ for $n \ge 3$ and some constant $C > 0$. Hence

$$\mathrm{Var}(\hat{\xi}_{\mathrm{tail}}) = o(n^{-1}). \tag{46}$$

We now turn to the bias. Note that

$$(E_\theta\hat{\xi}_{\mathrm{mid}} - \xi_{\mathrm{mid}})^2 \le 2(E_\theta\hat{\xi}_{\mathrm{mid1}} - \xi_{\mathrm{mid1}})^2 + 2(E_\theta\hat{\xi}_{\mathrm{mid2}} - \xi_{\mathrm{mid2}})^2$$

and

$$\begin{aligned}
(E_\theta\hat{\xi}_{\mathrm{mid1}} - \xi_{\mathrm{mid1}})^2 &\le E_\theta(\hat{\xi}_{\mathrm{mid1}} - \xi_{\mathrm{mid1}})^2 \\
&= E_\theta\biggl\{\biggl(\sum_{i=m_0+1}^{m_{k_*}} Y_i^2 - \lambda_{k_*}\biggr)_+ - \xi_{\mathrm{mid1}}\biggr\}^2 \\
&\le E_\theta\biggl(\sum_{i=m_0+1}^{m_{k_*}} Y_i^2 \\
&\qquad - \biggl(\frac{m_{k_*} - m_0}{n} + \xi_{\mathrm{mid1}}\biggr) - \biggl(\lambda_{k_*} - \frac{m_{k_*} - m_0}{n}\biggr)\biggr)^2 \\
&= \mathrm{Var}\biggl(\sum_{i=m_0+1}^{m_{k_*}} Y_i^2\biggr) + \biggl(\lambda_{k_*} - \frac{m_{k_*} - m_0}{n}\biggr)^2 \\
&\le \frac{4\xi_{\mathrm{mid1}}}{n} + \frac{5m_{k_*}\log m_{k_*}}{n^2}.
\end{aligned}$$

On the other hand, Lemma 1 shows that

$$(E_\theta\hat{\xi}_{\mathrm{mid2}} - \xi_{\mathrm{mid2}})^2 \le \biggl\{\sum_{i=m_{k_*}+1}^{m_J} \min\biggl(\frac{2\tau_{k_*,i}}{n}, \theta_i^2\biggr)\biggr\}^2,$$

and consequently the squared bias of the middle component satisfies

$$\begin{aligned}
(E_\theta\hat{\xi}_{\mathrm{mid}} - \xi_{\mathrm{mid}})^2 &\le \frac{8\xi_{\mathrm{mid1}}}{n} + \frac{10m_{k_*}\log m_{k_*}}{n^2} \\
&\quad + 2\biggl\{\sum_{i=m_{k_*}+1}^{m_J} \min\biggl(\frac{2\tau_{k_*,i}}{n}, \theta_i^2\biggr)\biggr\}^2.
\end{aligned} \tag{47}$$

Now consider the tail component. In this case Lemma 1 shows that the absolute bias of the tail component satisfies

$$|E_\theta\hat{\xi}_{\mathrm{tail}} - \xi_{\mathrm{tail}}| \le \sum_{i=m_J+1}^{\infty} \min\biggl(\frac{2\gamma_i\log n}{n}, \theta_i^2\biggr) + \sum_{i=m_J+1}^{\infty} \frac{4}{(2\pi\gamma_i)^{1/2}(\log n)^{1/2}n^{1+\gamma_i/2}}.$$

Note that $\gamma_i = 2(j+2)$ for $m_{J+j} + 1 \le i \le m_{J+j+1}$ and $j \ge 0$. Hence

$$\sum_{i=m_J+1}^{\infty} \frac{4}{(2\pi\gamma_i)^{1/2}(\log n)^{1/2}n^{1+\gamma_i/2}}$$

$$= \sum_{j=0}^{\infty} \sum_{i=m_{J+j}+1}^{m_{J+j+1}} \frac{4}{(4\pi(j+2))^{1/2}(\log n)^{1/2} n^{3+j}}$$

$$= \sum_{j=0}^{\infty} 2^j \frac{n^2}{(\log n)^2} \frac{4}{(4\pi(j+2))^{1/2}(\log n)^{1/2} n^{3+j}}$$

$$\le C n^{-1} (\log n)^{-5/2}$$

since $\sum_{j=0}^{\infty} \frac{2^j}{(j+2)^{1/2} n^j} \le \sum_{j=0}^{\infty} \frac{2^j}{(j+2)^{1/2} 3^j} < \infty$ whenever $n \ge 3$. Hence the squared bias of the tail component satisfies

$$(48)\quad (E_\theta \hat{\xi}_{\text{tail}} - \xi_{\text{tail}})^2 \le \left\{ \sum_{i=m_J+1}^{\infty} \min\left( \frac{2\gamma_i \log n}{n}, \theta_i^2 \right) + C n^{-1} (\log n)^{-5/2} \right\}^2.$$

We shall consider four separate cases.

*Case* 1. $p \ge 2$ and $\alpha > \frac{1}{4}$. In this case $k_* = 1$. Note that

$$\operatorname{Var}(\hat{Q}_{k_*}) = \operatorname{Var}(\hat{\xi}_0) + \operatorname{Var}(\hat{\xi}_{\text{mid}}) + \operatorname{Var}(\hat{\xi}_{\text{tail}}).$$

Note also that there exists a constant $C > 0$ such that for any $m \ge 1$

$$(49)\quad \sup_{\theta \in \Theta} \sum_{i=m}^{\infty} \theta_i^2 \le C m^{-2\alpha}.$$

It then follows from (45) that $\sup_{\theta\in\Theta} \operatorname{Var}(\hat{\xi}_{\text{mid}}) = o(n^{-1})$. This, together with (43) and (46), shows that $\sup_{\theta\in\Theta} V(\theta) = 4A(\Theta) n^{-1}(1 + o(1))$.

Now consider the bias. Note that $\hat{\xi}_0$ is an unbiased estimator of $\xi_0$. Hence

$$B^2(\theta) \le 2(E_\theta \hat{\xi}_{\text{mid}} - \xi_{\text{mid}})^2 + 2(E_\theta \hat{\xi}_{\text{tail}} - \xi_{\text{tail}})^2.$$

Note that in this case $\frac{n}{(\log n)^2} < m_{k_*} = 2m_0 \le \frac{2n}{(\log n)^2}$ and $\alpha > \frac{1}{4}$. Then (49) together with (47) and (48) yield that

$$\sup_{\theta\in\Theta} (E_\theta \hat{\xi}_{\text{mid}} - \xi_{\text{mid}})^2 \le \sup_{\theta\in\Theta} \left\{ \frac{8\xi_{\text{mid1}}}{n} + \frac{10 m_{k_*} \log m_{k_*}}{n^2} + 2\left( \sum_{i=m_{k_*}+1}^{m_J} \theta_i^2 \right)^2 \right\}$$
$$= o(n^{-1})$$

and

$$\sup_{\theta\in\Theta} (E_\theta \hat{\xi}_{\text{tail}} - \xi_{\text{tail}})^2 \le \sup_{\theta\in\Theta} \left\{ \sum_{i=m_J+1}^{\infty} \theta_i^2 + C n^{-1} (\log n)^{-5/2} \right\}^2 = o(n^{-1})$$

and consequently $\sup_{\theta\in\Theta} B^2(\theta) = o(n^{-1})$.

*Case* 2. $p \ge 2$ and $\alpha \le \frac{1}{4}$. In this case $m_{k_*}$ satisfies

$$\tfrac{1}{2} n^{2/(1+4\alpha)} (\log n)^{-1/(1+4\alpha)} < m_{k_*} \le n^{2/(1+4\alpha)} (\log n)^{-1/(1+4\alpha)}.$$

It then follows from (45) and (49) that

$$\sup_{\theta\in\Theta}\operatorname{Var}(\hat{\xi}_{\mathrm{mid}}) \le Cn^{-8\alpha/(1+4\alpha)}(\log n)^{-1/(1+4\alpha)}.$$

This together with (43) and (46) yield that $\sup_{\theta\in\Theta} V(\theta) \le Cn^{-8\alpha/(1+4\alpha)} \times (\log n)^{-1/(1+4\alpha)}$ for $\alpha < \frac{1}{4}$ and $\sup_{\theta\in\Theta} V(\theta) = 4A(\Theta)n^{-1}(1+o(1))$ for $\alpha = \frac{1}{4}$.

For the bias it follows from (47), (48) and (49) that

$$\begin{aligned}\sup_{\theta\in\Theta}(E_\theta\hat{\xi}_{\mathrm{mid}} - \xi_{\mathrm{mid}})^2 &\le \sup_{\theta\in\Theta}\left\{\frac{8\xi_{\mathrm{mid1}}}{n} + \frac{10m_{k_*}\log m_{k_*}}{n^2} + 2\left(\sum_{i=m_{k_*}+1}^{m_J}\theta_i^2\right)^2\right\}\\ &\le C\left(\frac{\log n}{n^2}\right)^{4\alpha/(1+4\alpha)},\end{aligned}$$

$$\sup_{\theta\in\Theta}(E_\theta\hat{\xi}_{\mathrm{tail}} - \xi_{\mathrm{tail}})^2 \le \sup_{\theta\in\Theta}\left\{\sum_{i=m_J+1}^{\infty}\theta_i^2 + Cn^{-1}(\log n)^{-5/2}\right\}^2 = o(n^{-8\alpha/(1+4\alpha)}),$$

and hence $\sup_{\theta\in\Theta} B^2(\theta) \le C(\frac{\log n}{n^2})^{4\alpha/(1+4\alpha)}$.

*Case* 3. $p < 2$ and $\alpha > \frac{1}{2p}$. This case is similar to Case 1. Note that in this case $k_* = 1$. Note also that for $p < 2$ the $L_p$ ball constraint (2) yields that for any $m \ge 1$

$$\sum_{i=m}^{\infty}\theta_i^2 \le M^2 m^{-2s}. \tag{50}$$

It then follows from (45) that $\sup_{\theta\in\Theta}\operatorname{Var}(\hat{\xi}_{\mathrm{mid}}) = o(n^{-1})$ and thus

$$\sup_{\theta\in\Theta} V(\theta) = 4A(\Theta)n^{-1}(1+o(1)).$$

We now turn to the bias. Note that it is straightforward to verify that for all $\theta \in L_p(\alpha, M)$ and all $j \ge 1$

$$\sum_{i=m_{k_*+j-1}+1}^{m_{k_*+j}}|\theta_i|^p \le M^p 2^{ps}2^{-jps}m_{k_*}^{-ps}. \tag{51}$$

Note also that for $m_{k_*+j-1}+1 \le i \le m_{k_*+j}$, $\tau_{k_*,i} = 2j$. Hence

$$\begin{aligned}\sum_{i=m_{k_*}+1}^{m_J}\min\left(\frac{2\tau_{k_*,i}}{n}, \theta_i^2\right) &= \sum_{j=1}^{J-k_*}\sum_{i=m_{k_*+j-1}+1}^{m_{k_*+j}}\min\left(\frac{4j}{n}, \theta_i^2\right)\\ &= \sum_{j=1}^{J-k_*}\frac{4j}{n}\sum_{i=m_{k_*+j-1}+1}^{m_{k_*+j}}\min\left(1, \theta_i^2\cdot\frac{n}{4j}\right)\end{aligned}$$

$$\le \sum_{j=1}^{J-k_*} \frac{4j}{n} \sum_{i=m_{k_*+j-1}+1}^{m_{k_*+j}} \min\left(1, \left\{\theta_i^2 \cdot \frac{n}{4j}\right\}^{p/2}\right),$$

where the last step follows from the facts $\min(1, \theta_i^2 \cdot \frac{n}{4j}) \le 1$ and $\frac{p}{2} \le 1$. Hence,

$$\sum_{i=m_{k_*}+1}^{m_J} \min\left(\frac{2\tau_{k_*,i}}{n}, \theta_i^2\right) \le \sum_{j=1}^{J-k_*} \left(\frac{4j}{n}\right)^{1-p/2} \sum_{i=m_{k_*+j-1}+1}^{m_{k_*+j}} |\theta_i|^p$$

$$\le \left\{M^p 2^{ps+2-p} \sum_{j=1}^{J-k_*} j^{1-p/2} 2^{-jps}\right\} \cdot m_{k_*}^{-ps} n^{p/2-1} \tag{52}$$

$$\le C m_{k_*}^{-ps} n^{-(1-p/2)}$$

for some constant $C > 0$, since $\sum_{j=1}^{\infty} j^{1-p/2} 2^{-jps} < \infty$. Similarly,

$$\sum_{i=m_J+1}^{\infty} \min\left(\frac{2\gamma_i \log n}{n}, \theta_i^2\right) \le C m_J^{-ps} n^{-(1-p/2)} (\log n)^{1-p/2}. \tag{53}$$

Note that $m_{k_*} = 2m_0 \ge n(\log n)^{-2}$ and $m_J \ge \frac{1}{4} n^2 (\log n)^{-2}$. Note also that in this case $\alpha p > \frac{1}{2}$. Hence $m_{k_*}^{-ps} n^{-(1-p/2)} = o(n^{-1/2})$ and $m_J^{-ps} n^{-(1-p/2)} \times (\log n)^{1-p/2} = o(n^{-1/2})$. Bounds in (52) and (53) together with (47) and (48) yield

$$\sup_{\theta\in\Theta} (E_\theta \hat{\xi}_{\text{mid}} - \xi_{\text{mid}})^2 = o(n^{-1}) \quad \text{and} \quad \sup_{\theta\in\Theta} (E_\theta \hat{\xi}_{\text{tail}} - \xi_{\text{tail}})^2 = o(n^{-1})$$

and consequently $\sup_{\theta\in\Theta} B^2(\theta) = o(n^{-1})$.

*Case* 4. $p < 2$ and $\alpha \le \frac{1}{2p}$. Note that in this case

$$\tfrac{1}{2} n^{p/(1+2ps)} (\log n)^{-1/(1+2ps)} < m_{k_*} \le n^{p/(1+2ps)} (\log n)^{-1/(1+2ps)}.$$

It follows from (45) and (49) that $\sup_{\theta\in\Theta} \operatorname{Var}(\hat{\xi}_{\text{mid}}) \le C n^{-(2-p/(1+2ps))} \times (\log n)^{-1/(1+2ps)}$. This together with (43) and (46) yield that $\sup_{\theta\in\Theta} V(\theta) \le C n^{-(2-p/(1+2ps))} (\log n)^{-1/(1+2ps)}$ for $\alpha < \frac{1}{2p}$ and $\sup_{\theta\in\Theta} V(\theta) = 4A(\Theta) n^{-1}(1 + o(1))$ for $\alpha = \frac{1}{2p}$.

On the other hand, (52) and (53) yield

$$\sum_{i=m_{k_*}+1}^{m_J} \min\left(\frac{2\tau_{k_*,i}}{n}, \theta_i^2\right) \le C n^{-1/2(2-p/(1+2ps))} (\log n)^{ps/(1+2ps)}, \tag{54}$$

$$\sum_{i=m_J+1}^{\infty} \min\left(\frac{2\gamma_i \log n}{n}, \theta_i^2\right) \le C n^{-p(\alpha+s)} (\log n)^{p(\alpha+s)}$$

$$= o(n^{-1/2(2-p/(1+2ps))}). \tag{55}$$

It now follows from (47) and (48) that

$$\sup_{\theta\in\Theta}(E_\theta\hat{\xi}_{\mathrm{mid}} - \xi_{\mathrm{mid}})^2 \le Cn^{-(2-p/(1+2ps)}(\log n)^{2ps/(1+2ps)},$$

$$\sup_{\theta\in\Theta}(E_\theta\hat{\xi}_{\mathrm{tail}} - \xi_{\mathrm{tail}})^2 = o(n^{-(2-p/(1+2ps))}),$$

and hence $\sup_{\theta\in\Theta} B^2(\theta) \le Cn^{-(2-p/(1+2ps))}(\log n)^{2ps/(1+2ps)}$.

REMARK. An inspection of the proof of Proposition 1 yields the following maximum mean squared error results for $\Theta = L_p(\alpha, M)$ and $\Theta = B^\alpha_{p,q}(M)$ which are useful for the proof of Theorem 2:

$$\sup_{\theta\in\Theta} E_\theta(\hat{\xi}_0 - \xi_0)^2 = 4A(\Theta)n^{-1}(1+o(1)), \tag{56}$$

$$\sup_{\theta\in\Theta} E_\theta(\hat{\xi}_{k_*} - \xi_{\mathrm{mid}})^2 = \begin{cases} o(n^{-1}), & \text{if } \alpha p_* > \frac{1}{2}, \\ Cn^{-(2-p_*/(1+2p_*s_*))}(\log n)^{2p_*s_*/(1+2p_*s_*)}, & \\ & \text{if } \alpha p_* \le \frac{1}{2}, \end{cases} \tag{57}$$

$$\sup_{\theta\in\Theta} E_\theta(\hat{\xi}_{\mathrm{tail}} - \xi_{\mathrm{tail}})^2 = \begin{cases} o(n^{-1}), & \text{if } \alpha p_* > \frac{1}{2}, \\ o(n^{-(2-p_*/(1+2p_*s_*))}), & \text{if } \alpha p_* \le \frac{1}{2}. \end{cases} \tag{58}$$

It should be stressed that in (57) $\hat{\xi}_{k_*}$ is the estimator defined by (24) which corresponds to a fixed Besov or $L_p$ ball.

4.3. *Proof of Theorem 2*. Let the estimator $\hat{Q}$ be given as in (31) and set $\Theta = L_p(\alpha, M)$. Note that $Q(\theta) = \xi_0 + \xi_{\mathrm{mid}} + \xi_{\mathrm{tail}}$. Note also that $\hat{\xi}_0$ is an unbiased estimate of $\xi_0$ and is independent of $\hat{\xi}_{\mathrm{mid}}$ and $\hat{\xi}_{\mathrm{tail}}$. Let the estimator $\hat{Q}$ be written as in (33). Then

$$\begin{aligned} E_\theta(\hat{Q} - Q(\theta))^2 &= E_\theta(\hat{\xi}_0 - \xi_0 + \hat{\xi}_{\mathrm{mid}} - \xi_{\mathrm{mid}} + \hat{\xi}_{\mathrm{tail}} - \xi_{\mathrm{tail}})^2 \\ &= E_\theta(\hat{\xi}_0 - \xi_0)^2 + E_\theta(\hat{\xi}_{\mathrm{mid}} - \xi_{\mathrm{mid}})^2 + E_\theta(\hat{\xi}_{\mathrm{tail}} - \xi_{\mathrm{tail}})^2 \\ &\quad + 2E_\theta(\hat{\xi}_{\mathrm{mid}} - \xi_{\mathrm{mid}})E_\theta(\hat{\xi}_{\mathrm{tail}} - \xi_{\mathrm{tail}}) \\ &\le E_\theta(\hat{\xi}_0 - \xi_0)^2 + 2E_\theta(\hat{\xi}_{\mathrm{mid}} - \xi_{\mathrm{mid}})^2 + 2E_\theta(\hat{\xi}_{\mathrm{tail}} - \xi_{\mathrm{tail}})^2. \end{aligned} \tag{59}$$

The difficulty lies in the analysis of $E_\theta(\hat{\xi}_{\mathrm{mid}} - \xi_{\mathrm{mid}})^2$ where $\hat{\xi}_{\mathrm{mid}}$ is given in (32), since the other terms $E_\theta(\hat{\xi}_0 - \xi_0)^2$ and $E_\theta(\xi_{\mathrm{tail}} - \xi_{\mathrm{tail}})^2$ satisfy (56) and (58). Set $\omega_k = \frac{6\sqrt{m_k \log n}}{n}$. A simple but important observation is that

$$\begin{aligned} (\hat{\xi}_{\mathrm{mid}} - \xi_{\mathrm{mid}})^2 &= \left(\max_{1\le k\le J}\{\hat{\xi}_k - \omega_k\} - \xi_{\mathrm{mid}}\right)^2 \\ &\le \min_{1\le k\le J}\{(\hat{\xi}_k - \omega_k - \xi_{\mathrm{mid}})^2\} + \sum_{k=1}^{J}[(\hat{\xi}_k - \omega_k - \xi_{\mathrm{mid}})_+]^2. \end{aligned}$$

Hence

$$
\begin{aligned}
E_\theta(\hat{\xi}_{\mathrm{mid}} - \xi_{\mathrm{mid}})^2 &\le \min_{1\le k\le J}\{E_\theta(\hat{\xi}_k - \omega_k - \xi_{\mathrm{mid}})^2\} \\
&\qquad + \sum_{k=1}^{J} E_\theta[(\hat{\xi}_k - \omega_k - \xi_{\mathrm{mid}})_+]^2 \\
&\le E_\theta(\hat{\xi}_{k_*} - \omega_{k_*} - \xi_{\mathrm{mid}})^2 \\
&\qquad + \sum_{k=1}^{J} E_\theta[(\hat{\xi}_k - \omega_k - \xi_{\mathrm{mid}})_+]^2
\end{aligned}
\tag{60}
$$

where $k_*$ is defined as in (23). The major difficulty in the analysis which follows is to show that the second term on the right-hand side of (60) is always negligible compared to the minimax risk. The first term is the dominant term and its analysis is made straightforward by the bounds given in (57).

For analysis of the second term on the right-hand side of (60), first consider the term $E_\theta[(\hat{\xi}_k - \omega_k - \xi_{\mathrm{mid}})_+]^2$. Let $\xi_{k,1} = \sum_{i=m_0+1}^{m_k} \theta_i^2$, $\xi_{k,2} = \sum_{i=m_k+1}^{m_J} \theta_i^2$, $\hat{\xi}_{k,1} = (\sum_{i=m_0+1}^{m_k} Y_i^2 - \lambda_k)_+$ and $\hat{\xi}_{k,2} = \sum_{i=m_k+1}^{m_J}[(Y_i^2 - \frac{\tau_{k,i}}{n})_+ - \mu_{k,i}]$. Note that it follows from the elementary inequality $(x+y)_+ \le (x)_+ + (y)_+$ for $x, y \in \mathbb{R}$ that

$$
\begin{aligned}
E_\theta[(\hat{\xi}_k - \omega_k - \xi_{\mathrm{mid}})_+]^2 &\le 2E_\theta[(\hat{\xi}_{k,1} - \xi_{k,1})_+]^2 \\
&\qquad + 2E_\theta[(\hat{\xi}_{k,2} - \omega_k - \xi_{k,2})_+]^2.
\end{aligned}
\tag{61}
$$

For the analysis of the first of these terms note that

$$
\begin{aligned}
[(\hat{\xi}_{k,1} - \xi_{k,1})_+]^2 &= \left[\left(\sum_{i=m_0+1}^{m_k} Y_i^2 - \lambda_k - \xi_{k,1}\right)_+\right]^2 \\
&= \left[\left(\frac{1}{n}\sum_{i=m_0+1}^{m_k} z_i^2 - \lambda_k + 2n^{-1/2}\sum_{i=m_0+1}^{m_k}\theta_i z_i\right)_+\right]^2.
\end{aligned}
$$

It then follows from the inequality

$$
[(x+y)_+]^2 \le 2[(x)_+]^2 + 2y^2 \tag{62}
$$

for any real numbers $x$ and $y$ that

$$
\begin{aligned}
E_\theta[(\hat{\xi}_{k,1} - \xi_{k,1})_+]^2 &\le 2E\left[\left(\frac{1}{n}\sum_{i=m_0+1}^{m_k} z_i^2 - \lambda_k\right)_+\right]^2 + 8n^{-1}E\left(\sum_{i=m_0+1}^{m_k}\theta_i z_i\right)^2 \\
&\le \frac{2}{n^2}E\left[\left(\sum_{i=m_0+1}^{m_k} z_i^2 - n\lambda_k\right)_+\right]^2 + 8\xi_{k,1}n^{-1}.
\end{aligned}
$$

Set $m = m_k - m_0$. It then follows from Theorem 2.1 of [21] that

$$E\left[\left(\sum_{i=m_0+1}^{m_k} z_i^2 - n\lambda_k\right)_+\right]^2 \le 8\left(\frac{2\sqrt{m\log m}+m}{2\sqrt{m\log m}+2}\right)^2 P(X_m \ge m + 2\sqrt{m\log m}\,),$$

where $X_m$ is a central chi-square random variable with $m$ degrees of freedom. It then follows from Lemma 2 of [8] on the tail probability bounds of the chi-square distribution and by noting $\log(1+x) \le x - \frac{1}{2}x^2 + \frac{1}{3}x^3$ for all $x \ge 0$ that

$$\begin{aligned} P(X_m \ge m + 2\sqrt{m\log m}\,) &\le \frac{1}{2}\exp\left(-\frac{m}{2}\left[2\sqrt{\frac{\log m}{m}} - \log\left(1 + 2\sqrt{\frac{\log m}{m}}\right)\right]\right) \\ &\le \frac{1}{2}\exp\left(-\log m + \frac{4(\log m)^{3/2}}{3m^{1/2}}\right) \le \frac{3}{m}, \end{aligned}$$

where the last inequality follows from the fact that $\frac{(\log m)^{3/2}}{m^{1/2}}$ attains its maximum at $m = e^3$ and hence $\frac{1}{2}\exp(\frac{4(\log m)^{3/2}}{3m^{1/2}}) \le 3$. Therefore

$$E\left[\left(\sum_{i=m_0+1}^{m_k} z_i^2 - n\lambda_k\right)_+\right]^2 \le 8\left(\frac{2\sqrt{m\log m}+m}{2\sqrt{m\log m}+2}\right)^2 \cdot \frac{3}{m} \le \frac{24}{\log m}$$

and hence

$$E_\theta[(\hat{\xi}_{k,1} - \xi_{k,1})_+]^2 \le \frac{48}{n^2\log(m_k - m_0)} + \frac{8\xi_{k,1}}{n}. \tag{63}$$

We now turn to $E_\theta[(\hat{\xi}_{k,2} - \omega_k - \xi_{k,2})_+]^2$ and show that the follow bound holds.

Lemma 4. *For some constant $C > 0$ and for all sufficiently large $n$*

$$E_\theta[(\hat{\xi}_{k,2} - \omega_k - \xi_{k,2})_+]^2 \le Cn^{-1}(\log n)^{-5/4} + 4\xi_{k,2}n^{-1}\log_2 n. \tag{64}$$

Proof. Set $D_{k,2} = E_\theta[(\hat{\xi}_{k,2} - \omega_k - \xi_{k,2})_+]^2$. Note that

$$\hat{\xi}_{k,2} = \sum_{i=m_k+1}^{m_J}\left[\left(Y_i^2 - \frac{\tau_{k,i}}{n}\right)_+ - \mu_{k,i}\right] = \sum_{j=1}^{J-k}\sum_{i=m_{j+k-1}+1}^{m_{j+k}}\left[\left(Y_i^2 - \frac{2j}{n}\right)_+ - \mu_{k,i}\right].$$

Set $\eta_j = \sum_{i=m_{j+k-1}+1}^{m_{j+k}} \theta_i^2$ for $1 \le j \le J-k$. It follows from Lemma 2 that for a fixed value of $\eta_j$ on a block $m_{j+k-1}+1 \le i \le m_{j+k}$, $\sum_{i=m_{j+k-1}+1}^{m_{j+k}}(Y_i^2 - \frac{2j}{n})_+$ is stochastically maximized when $\theta_{m_{j+k-1}+1} = \eta_j^{1/2}$ and the remaining $\theta_i = 0$.

Hence

$$D_{k,2} \le E\Bigg\{\Bigg(\sum_{j=1}^{J-k}\Bigg[\Bigg(\eta_j + 2\eta_j^{1/2} n^{-1/2} z_{m_{j+k-1}+1} + \frac{1}{n} z_{m_{j+k-1}+1}^2 - \frac{2j}{n}\Bigg)_+ - \mu_{k,m_{j+k-1}+1}\Bigg] + \sum_{j=1}^{J-k}\sum_{i=m_{j+k-1}+2}^{m_{j+k}}\Bigg[\Bigg(\frac{1}{n}z_i^2 - \frac{2j}{n}\Bigg)_+ - \mu_{k,i}\Bigg] - \omega_k - \xi_{k,2}\Bigg)_+\Bigg\}^2.$$

Noting $\sum_{j=1}^{J-k}\eta_j = \xi_{k,2}$, it then follows from the fact that $(x+y)_+ \le (x)_+ + (y)_+$ and (62) that

$$\begin{aligned} D_{k,2} &\le E\Bigg\{\Bigg(\sum_{j=1}^{J-k} 2\eta_j^{1/2} n^{-1/2}(z_{m_{j+k-1}+1})_+ + \sum_{j=1}^{J-k}\sum_{i=m_{j+k-1}+1}^{m_{j+k}}\Bigg[\Bigg(\frac{1}{n}z_i^2 - \frac{2j}{n}\Bigg)_+ - \mu_{k,i}\Bigg] - \omega_k\Bigg)_+\Bigg\}^2 \\ &\le 2n^{-2}E\Bigg\{\Bigg(\sum_{j=1}^{J-k}\sum_{i=m_{j+k-1}+1}^{m_{j+k}}[(z_i^2 - 2j)_+ - n\mu_{k,i}] - n\omega_k\Bigg)_+\Bigg\}^2 \\ &\quad + 2n^{-1}E\Bigg\{\sum_{j=1}^{J-k} 2\eta_j^{1/2}(z_{m_{j+k-1}+1})_+\Bigg\}^2. \end{aligned} \tag{65}$$

It is easy to see that the second term

$$\begin{aligned} 2n^{-1}E\Bigg\{\sum_{j=1}^{J-k} 2\eta_j^{1/2}(z_{m_{j+k-1}+1})_+\Bigg\}^2 &\le 2n^{-1}J\sum_{j=1}^{J-k} 4\eta_j E\{(z_{m_{j+k-1}+1})_+\}^2 \\ &= 4J\xi_{k,2}n^{-1}. \end{aligned} \tag{66}$$

We now use Lemmas 1 and 3 to bound $E\{(\sum_{j=1}^{J-k}\sum_{i=m_{j+k-1}+1}^{m_{j+k}}[(z_i^2 - 2j)_+ - n\mu_{k,i}] - n\omega_k)_+\}^2$. First note that equation (37) in Lemma 1 yields that

$$\mathrm{Var}((z_i^2 - 2j)_+) \le [16(2j)^{-1/2} - 9(2j)^{-3/2} + 9(2j)^{-5/2}] \cdot (2\pi)^{-1/2} e^{-j}$$

and consequently

$$\begin{aligned} V_n &\equiv \sum_{j=1}^{J-k}\sum_{i=m_{j+k-1}+1}^{m_{j+k}} (\mathrm{Var}(z_i^2 - 2j)_+) \\ &\le m_k \sum_{j=1}^{J-k}[16(2j)^{-1/2} - 9(2j)^{-3/2} + 9(2j)^{-5/2}] \cdot (2\pi)^{-1/2}\Bigg(\frac{e}{2}\Bigg)^{-j} \le 9m_k, \end{aligned}$$

where the last step follows from the fact that

$$\sum_{j=1}^{\infty}[16(2j)^{-1/2} - 9(2j)^{-3/2} + 9(2j)^{-5/2}] \cdot (2\pi)^{-1/2}\left(\frac{e}{2}\right)^{-j} < 9,$$

which can be verified by direct calculations. It then follows from Lemma 3 with $\gamma = 4$ that for all sufficiently large $n$

$$E\left\{\left(\sum_{j=1}^{J-k}\sum_{i=m_{j+k-1}+1}^{m_{j+k}}[(z_i^2-2j)_+ - n\mu_{k,i}] - n\omega_k\right)_+\right\}^2 \le Cm_k(\log n)^{3/4}n^{-1}, \tag{67}$$

where $C > 0$ is a constant. Noting $m_k \le m_J \le \frac{n^2}{(\log n)^2}$ and $J \le \log_2 n$, (65), (66) and (67) together yield that

$$\begin{aligned} E_\theta[(\hat{\xi}_{k,2} - \omega_k - \xi_{k,2})_+]^2 &\le Cm_k(\log n)^{3/4}n^{-3} + 4J\xi_{k,2}n^{-1} \\ &\le Cn^{-1}(\log n)^{-5/4} + 4\xi_{k,2}n^{-1}\log_2 n \end{aligned}$$

and Lemma 4 is thus proved. □

We now return to the proof of Theorem 2. Lemma 4 together with (61) and (63) yield that for all sufficiently large $n$

$$\begin{aligned} \sum_{k=1}^{J} E_\theta[(\hat{\xi}_k - \omega_k - \xi_{\mathrm{mid}})_+]^2 &\le CJn^{-2}(\log n)^{-1} + CJn^{-1}(\log n)^{-5/4} \\ &\quad + 4J\xi_{\mathrm{mid}}n^{-1}\log_2 n \\ &\le C\{n^{-2} + n^{-1}(\log n)^{-1/4} + \xi_{\mathrm{mid}}n^{-1}(\log n)^2\}. \end{aligned}$$

It then follows from (49) for $p \ge 2$ and (50) for $p < 2$ that

$$\sup_{\theta\in\Theta}\sum_{k=1}^{J} E_\theta[(\hat{\xi}_k - \omega_k - \xi_{\mathrm{mid}})_+]^2 = o(n^{-1}) \tag{68}$$

and is thus negligible relative to the minimax risk.

The rest of the proof is now straightforward. It follows from (59) and (60) that

$$\begin{aligned} E_\theta(\hat{Q} - Q(\theta))^2 &\le E_\theta(\hat{\xi}_0 - \xi_0)^2 + 2E_\theta(\hat{\xi}_{k_*} - \omega_{k_*} - \xi_{\mathrm{mid}})^2 \\ &\quad + 2\sum_{k=1}^{J} E_\theta[(\hat{\xi}_k - \omega_k - \xi_{\mathrm{mid}})_+]^2 + 2E_\theta(\hat{\xi}_{\mathrm{tail}} - \xi_{\mathrm{tail}})^2 \end{aligned} \tag{69}$$

$$\leq \{E_\theta(\hat{\xi}_0 - \xi_0)^2 + 4E_\theta(\hat{\xi}_{k_*} - \xi_{\text{mid}})^2 + 2E_\theta(\hat{\xi}_{\text{tail}} - \xi_{\text{tail}})^2\}$$

$$+ 4\omega_{k_*}^2 + 2\sum_{k=1}^{J} E_\theta[(\hat{\xi}_k - \omega_k - \xi_{\text{mid}})_+]^2.$$

The remainder of the proof can be separated into two cases, $\alpha p_* \leq \frac{1}{2}$ and $\alpha p_* > \frac{1}{2}$. First consider the case when $\alpha p_* \leq \frac{1}{2}$. In this case it follows from the definition of $m_{k_*}$ given in (23) that

$$\omega_{k_*}^2 = \frac{36 m_{k_*} \log n}{n^2} \leq 72 n^{-(2-p_*/(1+2p_* s_*))} (\log n)^{2p_* s_*/(1+2p_* s_*)}. \tag{70}$$

For this case the theorem now immediately follows from (69), (68), (70) and (56)–(58).

For the case $\alpha p_* > \frac{1}{2}$ first note that $k_* = 1$ and $\omega_{k_*} = 2m_0 = o(n^{-1})$. The theorem then immediately follows in this case from this observation, (69), (68) and (56)–(58).

4.4. *Proof of Theorem 1.* We divide the proof into two cases, $p \geq 2$ and $0 < p < 2$, which correspond to the cases $p_* = 2$ and $p_* < 2$, respectively. The case where $p \geq 2$ is standard but we include a brief outline for the sake of completeness. In this case we apply Theorem 2.1 of [14] combined with Theorem 4 of [11]. Let

$$\omega(\delta) = \sup\left\{Q(\theta) : \sum_{i=1}^{\infty} \theta_i^2 \leq \delta^2, \theta \in L_p(\alpha, M)\right\}$$

be the modulus of continuity introduced in [13]. For small $\delta$ let $N \sim \delta^{-2/(4\alpha+1)}$. Let $\theta = (\theta_1, \theta_2, \ldots)$, where $\theta_i = c\delta^{(2\alpha+1)/(4\alpha+1)}$ for $i = 1, \ldots, n$ and otherwise $\theta_i = 0$. It is easy to check that $\theta \in L_p(\alpha, M)$ for sufficiently small $c > 0$. Simple calculations then show that $\omega(\delta) \geq D\delta^{4\alpha/(4\alpha+1)}$ for some $D > 0$.

It then follows from Theorem 2.1 of [14] and Theorem 4 of [11] that if $\hat{Q}$ satisfies (10) then

$$\sup_{\theta \in \Theta} (E_\theta \hat{Q} - Q(\theta))^2 \geq \omega^2\left(d\frac{\sqrt{\log n}}{n}\right) - 2\sqrt{C} n^{-\gamma/2} \omega\left(d\frac{\sqrt{\log n}}{n}\right) n^{d^2}. \tag{71}$$

Equation (11) follows by taking a sufficiently small $d$.

We now turn to the case where $p < 2$ and $\alpha < \frac{1}{2p}$, in which case $p_* = p$ and $s_* = s$. The proof follows a similar argument to one given in [10] for minimax lower bounds. The main idea is to place a prior on the union of the zero vector and the vertices of a suitable collection of hypercubes. The constrained risk inequality given in Theorem 4 of [11] can then be used to yield a lower bound for the maximum mean squared error over the vertices of the hypercubes given an upper bound on the mean squared error at the origin.

More precisely, let $\Theta_{k,m}$ be the union of the zero vector $\theta_0 = (0, 0, \ldots)$ and the collection of vectors which have exactly $k$ nonzero coordinates equal to $\frac{1}{\sqrt{n}}$ in the first $m$ coordinates and are otherwise equal to zero. It is straightforward to check that $\Theta(k,m) \subset L_p(\alpha, M)$ when $m = n^{p/(1+2ps)}(\log n)^{-1/(1+2ps)}$ and $k = \sqrt{bm \log m}$ for sufficiently small constant $b > 0$.

As in [10] let $\mathcal{I}(k,m)$ be the class of all subsets of $\{1, \ldots, m\}$ of $k$ elements and for $I \in \mathcal{I}(k,m)$ let $\theta_I \in \Theta_{k,m}$ be the vector where the $j$th coordinate is zero if $j \notin I$ and is equal to $\frac{1}{\sqrt{n}}$ for $j \in I$.

Let $\psi_\mu$ be the density of a normal distribution with mean $\mu$ and variance $\frac{1}{n}$. And for $I \in \mathcal{I}(k,m)$ let $g_I(y_1, \ldots, y_m) = \prod_{j=1}^m \psi_{\mu_j}(y_j)$, where $\mu_j = \frac{1}{\sqrt{n}} \mathbb{1}(j \in I)$. Finally let $g = \frac{1}{\binom{m}{k}} \sum_{I \in \mathcal{I}(k,m)} g_I$ and $f$ be the density of $m$ independent normal random variables each with mean 0 and variance $\frac{1}{n}$. Note that a similar mixture prior was used in [2] to give lower bounds in a nonparametric testing problem.

The application of the constrained risk inequality of [11] requires an upper bound on the chi-squared distance between $f$ and $g$. Cai and Low [10] shows that $\int \frac{g^2}{f} = Ee^J$ where $J$ has the hypergeometric distribution $P(J = j) = \frac{\binom{k}{j}\binom{m-k}{k-j}}{\binom{m}{k}}$. Now note from [16], page 59, that

$$P(J = j) \le \binom{k}{j} \left(\frac{k}{m}\right)^j \left(1 - \frac{k}{m}\right)^{k-j} \left(1 - \frac{k}{m}\right)^{-k}.$$

For $k = \sqrt{bm \log m}$, $(1 - \frac{k}{m})^{-k} \le e^{k^2/m} \le m^b$ and it follows that

$$\int \frac{g^2}{f} \le m\left(1 + (e-1)\frac{k}{m}\right)^k \le m e^{(e-1)k^2/m} \le m^{be}. \tag{72}$$

The constrained risk inequality in Theorem 4 of [11] then yields that if for any $\varepsilon > 0$, $E_f(\hat{Q} - Q(\theta_0))^2 \le \frac{m^{1-\varepsilon}}{n^2}$, then for $b \le \frac{\varepsilon}{e}$

$$\begin{aligned} \left(E_g \hat{Q} - k\frac{\rho^2}{n}\right)^2 &\ge \frac{k^2}{n^2} - 2m^{be/2}\frac{k}{n}\frac{m^{(1-\varepsilon)/2}}{n} = \frac{k^2}{n^2}(1 + o(1)) \\ &= \frac{bm \log m}{n^2}(1 + o(1)). \end{aligned} \tag{73}$$

Hence for some constant $C_1 > 0$

$$\begin{aligned} &\inf_{\hat{Q}} \sup_{\theta \in L_p(\alpha, M)} (E_\theta \hat{Q} - Q(\theta))^2 \\ &\qquad \ge \inf_{\hat{Q}} \sup_{\theta \in \Theta_{k,m}} (E_\theta \hat{Q} - Q(\theta))^2 \ge \frac{bm \log m}{n^2} \end{aligned}$$

$$\geq C_1 n^{p/(1+2ps)-2} (\log n)^{2ps/(1+2ps)} (1 + o(1)).$$

## REFERENCES

[1] Baraud, Y. (2000). Model selection for regression on a fixed design. *Probab. Theory Related Fields* **117** 467–493. MR1777129

[2] Baraud, Y. (2002). Non-asymptotic minimax rates of testing in signal detection. *Bernoulli* **8** 577–606. MR1935648

[3] Barron, A. R., Birgé, L. and Massart, P. (1999). Risk bound for model selection via penalization. *Probab. Theory Related Fields* **113** 301–413. MR1679028

[4] Bickel, P. J., Klaassen, C. A. J., Ritov, Y. and Wellner, J. A. (1993). *Efficient and Adaptive Estimation for Semiparametric Models*. Johns Hopkins Univ. Press, Baltimore. MR1245941

[5] Bickel, P. J. and Ritov, Y. (1988). Estimating integrated squared density derivatives: Sharp best order of convergence estimates. *Sankhyā Ser. A* **50** 381–393. MR1065550

[6] Birgé, L. and Massart, P. (1997). From model selection to adaptive estimation. In *Festschrift for Lucien Le Cam*: *Research Papers in Probability and Statistics* (D. Pollard, E. Torgersen and G. Yang, eds.) 55–87. Springer, New York. MR1462939

[7] Brown, L. D. (1986). *Fundamentals of Statistical Exponential Families with Applications in Statistical Decision Theory*. IMS, Hayward, CA. MR0882001

[8] Cai, T. (1999). Adaptive wavelet estimation: A block thresholding and oracle inequality approach. *Ann. Statist.* **27** 898–924. MR1724035

[9] Cai, T. and Low, M. (2005). On adaptive estimation of linear functionals. *Ann. Statist.* **33** 2311–2343. MR2211088

[10] Cai, T. and Low, M. (2005). Nonquadratic estimators of a quadratic functional. *Ann. Statist.* **33** 2930–2956.

[11] Cai, T., Low, M. and Zhao, L. (2006). Tradeoffs between global and local risks in nonparametric function estimation. *Bernoulli*. To appear.

[12] Donoho, D. L. and Liu, R. C. (1991). Geometrizing rates of convergence. II, III. *Ann. Statist.* **19** 633–667, 668–701. MR1105839

[13] Donoho, D. L. and Nussbaum, M. (1990). Minimax quadratic estimation of a quadratic functional. *J. Complexity* **6** 290–323. MR1081043

[14] Efromovich, S. Y. and Low, M. (1996). On optimal adaptive estimation of a quadratic functional. *Ann. Statist.* **24** 1106–1125. MR1401840

[15] Fan, J. (1991). On the estimation of quadratic functionals. *Ann. Statist.* **19** 1273–1294. MR1126325

[16] Feller, W. (1968). *An Introduction to Probability Theory and Its Applications* **1**, 3rd ed. Wiley, New York. MR0228020

[17] Feller, W. (1971). *An Introduction to Probability Theory and Its Applications* **2**, 2nd ed. Wiley, New York. MR0270403

[18] Gayraud, G. and Tribouley, K. (1999). Wavelet methods to estimate an integrated quadratic functional: Adaptivity and asymptotic law. *Statist. Probab. Lett.* **44** 109–122. MR1706448

[19] Hall, P. and Heyde, C. C. (1980). *Martingale Limit Theory and Its Application*. Academic Press, New York. MR0624435

[20] Hardy, G., Littlewood, J. E. and Pólya, G. (1952). *Inequalities*, 2nd ed. Cambridge Univ. Press. MR0046395

[21] Johnstone, I. M. (2001). Chi-square oracle inequalities. In *State of the Art in Probability and Statistics. Festschrift for Willem R. van Zwet* (M. de Gunst,

C. Klaassen and A. van der Waart, eds.) 399–418. IMS, Beachwood, OH. MR1836572

[22] Johnstone, I. M. (2002). Function estimation and Gaussian sequence models. Draft of a monograph. Available at stat.stanford.edu/~imj/baseb.pdf.

[23] Klemelä, J. (2006). Sharp adaptive estimation of quadratic functionals. *Probab. Theory Related Fields* **134** 539–564. MR2214904

[24] Klemelä, J. and Tsybakov, A. B. (2001). Sharp adaptive estimation of linear functionals. *Ann. Statist.* **29** 1567–1600. MR1891739

[25] Laurent, B. and Massart, P. (2000). Adaptive estimation of a quadratic functional by model selection. *Ann. Statist.* **28** 1302–1338. MR1805785

[26] Lepski, O. V. (1990). On a problem of adaptive estimation in Gaussian white noise. *Theory Probab. Appl.* **35** 454–466. MR1091202

[27] Lepski, O. V. and Spokoiny, V. G. (1997). Optimal pointwise adaptive methods in nonparametric estimation. *Ann. Statist.* **25** 2512–2546. MR1604408

[28] Mardia, K. V. (1972). *Statistics of Directional Data*. Academic Press, London. MR0336854

[29] Tribouley, K. (2000). Adaptive estimation of integrated functionals. *Math. Methods Statist.* **9** 19–38. MR1772223

[30] Watson, G. S. (1983). *Statistics on Spheres*. Wiley, New York. MR0709262

Department of Statistics
The Wharton School
University of Pennsylvania
Philadelphia, Pennsylvania 19104-6340
USA
E-mail: tcai@wharton.upenn.edu
lowm@wharton.upenn.edu